\documentclass[11pt, leqno]{amsart}

\usepackage{amsmath}
\usepackage{amssymb}

\begin{document}

\title{ On three dimensional conformally flat almost cosymplectic
          manifolds}
\author{Piotr Dacko}


\begin{abstract}
In the paper there are described new examples of conformally flat
three dimensional almost cosymplectic manifolds.  All these manifolds
form a class which was completely characterized.      
\end{abstract}

\maketitle

\newtheorem{proposition}{Proposition}
\newtheorem{theorem}{Theorem}

\section{ Introduction}

It is a difficult problem to construct explicit examples of of almost cosymplectic manifolds endowed with  
Riemannian metrics which satisfy some  classical curvature conditions like e.g.  the Einstein condition, conformall flatness, locall symetricity etc.  Known  explicit examples are of a very special structure.  Many of them are some Lie groups endowed with a left invariant almost cosymplectic structures \cite{CLM-2}, \cite{CFL}, \cite{DO2}.  However with no evident connection to the  mentioned above metric properties.  In this view any new explicit examples are of no doubt of much worth as they may give a new impulse for the further development of the hole theory.

For an almost cosymplectic manifolds there are non-existence theorems.  The oldest result of this kind is the theorem due to Z.Olszak \cite{O1}, \cite{O2}  which asserts that the curvature of a constant curvature almost cosymplectic manifold is zero and manifold is cosymplectic.  For the almost cosymplectic manifold with K\"ahler leaves (see next section) holds even stronger result: in dimension $\geq 5$ there are no conformally flat metric \cite{DO}.  The case of the dimension $3$ is different.  
We remind here the first explicit example of a  conformally flat almost cosymplectic three manifold given in \cite{DO}: $p=(x,y,z)\in U \subset \mathbb R^3$,
\begin{eqnarray}
\label{z^2}
\nonumber & g = z^2dx^2 + \dfrac{e^{2ax}}{z^2}dy^2 + dz^2, & \\
& \xi = \dfrac{\partial}{\partial z}, \quad\quad \eta = dz, & \\ [+4pt]
\nonumber &  \varphi \dfrac{\partial}{\partial x} = 
\dfrac{z^2}{e^{ax}}\dfrac{\partial}{\partial y}, \quad\quad
\varphi \dfrac{\partial}{\partial y} = 
- \dfrac{e^{ax}}{z^2} \dfrac{\partial}{\partial x}. &
\end{eqnarray}

Remarkable the almost contact structure in this example is non-normal therefore this it is  a non-cosymplectic conformally flat manifold.  This is important as it is easy to construct an example of conformally flat {\em cosymplectic} three manifold.  Even more to classify locally all such manifolds (cf. sect. 4).

There are no explicit examples in dimensions $\geq 5$ and it is still the open problem whether there exist conformally flat non-flat almost cosymplectic manifolds.  In \cite{O1} Z. Olszak proved that the scalar  curvature $s$ of conformally flat almost cosymplectic manifold of dimension $\geq 5$ is non-positive which was later improved by H. Endo \cite{E} who found the following inequality  
\begin{equation*}
\rho(X) + \rho(\varphi X) \leq -\dfrac{s}{n(2n-3)},
\end{equation*}
for the Ricci curvature $\rho$ of arbitrary unit vector $X\perp \xi$.

\section{Preliminaries}
An almost contact metric structure consists of four tensor fields customary denoted by $\varphi$, $\xi$, $\eta$ and $g$ where $\varphi$ is a $(1,1)$ tensor field, $\xi$ is a vector field, $\eta$ a $1$-form and $g$ is a Riemannian metric.  Moreover one requires the following relations must be satisfied
\begin{equation}
\label{compat}
\begin{array}{c}
 \varphi^2 = -Id + \eta\otimes\xi, \quad\quad  \eta(\xi) = 1,  \\ [+4pt]
 g(\varphi X,\varphi Y) = g(X,Y) - \eta(X)\eta(X), 
\end{array}
\end{equation}
where $X,Y$ are arbitrary vector fields.  

To each almost contact metric manifold $M$, i.e.  manifold endowed with an almost contact metric structure, is associated  a $2$-form $\varPhi(X,Y)= g(\varphi X,Y)$ usually called a fundamental form.  We mention  that our definition of $\varPhi$ differs by sign of that given in \cite{B} where $\varPhi(X,Y)=g(X,\varphi Y)$.  However this is explained by historical reasons.  
  
An almost contact metric ma\-ni\-fold $M$ is always odd-di\-men\-sio\-nal, $\dim M$ 
$= 2n+1$.  The fundamental form $\varPhi$ is of maximal possible rank $2n$ as its kernel, at each point,  consists of vectors  $c\xi$, $c= const$.  Moreover $\omega= \eta\wedge\varPhi^n$ is a non-vanishing everywhere $(2n+1)$-form hence $M$ is orientable.

Considering the behavior of the exterior differentials of the forms $\eta$ and $\varPhi$ we fall through different classes of almost contact metric manifolds.  The more exhaustive studied class of manifolds nowadays  are {\em contact metric} manifolds.  The terminology is  explained
 by the fact that $\eta$ is a contact form $\eta\wedge (d\eta)^n\neq 0$.  The basic reference to this theory is a monograph \cite{B}.  In a some sense opposite direction we have {\em almost cosymplectic} manifolds as they were defined by the conditions that both $\eta$ and $\varPhi$ are closed \cite{GY}.  

Let $M$ be an almost cosymplectic manifold.  Then we have the following fundamental identity \cite{O1}
\begin{equation*}
  (\nabla_{\varphi X}\varphi)\varphi Y
 + (\nabla_X\varphi)Y -\eta(Y)\nabla_{\varphi X}\xi = 0.
\end{equation*}
which implies that $\nabla_\xi\varphi = 0$ and $ \nabla_\xi\xi = 0$
hence any integral curve of $\xi$ is a geodesic.  

Let denote by $\mathcal F$ a foliations of $M$ defined by $\mathcal D = \ker\eta$.  We fix a leaf $N\in \mathcal F$.  A form 
\begin{equation*}
\Omega = i^*\varPhi,
\end{equation*} 
$i$ being an inclusion map, is a symplectic form on $N$ so $(N,\Omega)$ is a symplectic manifold.  Even more $(N,\Omega)$ can be endowed with an inherited almost Hermitian structure $(J,G)$ in the way that $\Omega$ becomes a fundamental form of this almost Hermitian structure.  Therefore $(N, J, G)$ may be considered as an almost K\"ahler manifold. 

 We set by definition
\begin{equation}
\label{Jdef}
  i_*(J\bar X) = \varphi i_*(\bar X),
\end{equation}
for any tangent to $N$ vector field.  We have  $i_*(\bar X) = \mathcal D_{i(p)}$ as $N$ is an integral submanifold of $\mathcal D$.  Note
\begin{equation*}
   \varphi^2 i_*(\bar X) = -i_*(\bar X) + \eta(i_*(\bar X))\xi = -i_*(\bar X).
\end{equation*}
  Thus $\varphi^2 |_{\mathcal D_p}= -Id |_{\mathcal D_p}$ and the linear algebra arguments imply $\varphi(\mathcal D_p) = \mathcal D_p$.  The identity above also shows that $J$ is an almost complex structure on $N$.  Now for the metric $G$ we set $G=i^*g$.  Therefore $(N,G)$ is a Riemannian hypersurface in $M$.  Let $\bar X, \bar Y$ be arbitrary tangent to $N$ vector fields.  Then
\begin{eqnarray*}
  G(J\bar X, J\bar Y) & = & g(i_*(J\bar X), i_*(J\bar Y)) = g(\varphi i_*(\bar X), \varphi i_*(\bar Y)) \\
  & = & g( i_*(\bar X), i_*(\bar Y)) = G(\bar X, \bar Y),
\end{eqnarray*}
where we have used (\ref{Jdef}), the definition of $G$ and (\ref{compat}).   In similar way we can show that
\begin{equation*}
\Omega(\bar X, \bar Y) = G(J\bar X, \bar Y).
\end{equation*}
Summing up all above we have
\begin{eqnarray*}
& J^2 = -Id, \quad\quad G(JX,JY) = G(X, Y), & \\
& \Omega(X,Y) = G(JX,Y), \quad\quad d\Omega = 0,
\end{eqnarray*}
so $(J,G)$ is an almost K\"ahler structure on $N$ with $\Omega$ as fundamental form.  Of course our construction do not depend on the choice of  a leaf in the sense that each leaf can be endowed with an almost K\"ahler structure in the way described above.  However in general these structures on different leaves are different.  
In the case that these induced structures on any leaf are K\"ahlerian  the manifold $M$ is said to be almost cosymplectic with K\"ahlerian leaves \cite{O2}.  Note that each three dimensional almost cosymplectic manifold clearly has K\"ahlerian leaves as the leaves are two dimensional.

Let $\nabla$ and $\bar \nabla$ denotes the Levi-Civitta connections on $M$ and $N$ resp.  then we have the usual Gauss decomposition formula
\begin{equation*}
  \nabla_{\bar X}\bar Y = \bar \nabla_{\bar X}\bar Y + h(\bar X, \bar Y)\bar n,
\end{equation*}
where $h$ is a second fundamental form of the Riemannian hypersurface $N$ and $\bar n$ stands for the normal vector field.  The field will be determined uniquely if we choose orientation on both $N$ and $M$ requiring that for a given positively oriented frame $(\bar X_1,\ldots,\bar X_{2n})$ of $T_pN$ the frame $(\bar n, \bar X_1,\ldots, \bar X_{2n})$ of $T_pM$ is also positively oriented.  Quite naturally we take on $M$ an orientation given by the equivalence class of the form $\omega= \eta\wedge\varPhi^n$  and for $N$  those one determined by the almost complex structure $J$.   With these assumptions we have $\bar n = \xi|_N$ as $\xi$ is unit vector field and everywhere orthogonal to any tangent space of $N$.

For the Weingarten operator of $N$ we have
\begin{equation*}
 S\bar X = - \nabla_{\bar X}\bar n = -\nabla_{\bar X}\xi.
\end{equation*}
   Again we see that formula for Weingarten operators of different leaves is exactly the same.   This observation  suggests to introduce a  tensor field  $A$ as follows
\begin{equation*}
AX = -\nabla_X\xi.
\end{equation*}

We note some properties of $A$.
On the leaf $N$ the tensor field $A$ and the operator $S$ are related by 
\begin{equation}
\label{as}  A i_*(\bar X) = i_*(S\bar X),
\end{equation}
$A$ is symmetric $g(AX,Y) = g(X, AY)$ and anti-commutes with $\varphi$ 
\begin{equation}
\label{aphi}  \varphi A + A\varphi = 0.  
\end{equation}
As the vector field $\xi$ is geodesic we have
$A\xi = -\nabla_\xi\xi =0$ which implies $\eta(AX) =0 $.  From (\ref{aphi}) it follows that
\begin{equation*}
\label{tra}
 A\varphi X = -\lambda \varphi X \quad \mbox{if} \quad AX = \lambda X.
\end{equation*}
Therefore the spectrum of $A$ always is of the form 
\begin{equation*}
(0,\lambda_1,\lambda_2,\ldots,\lambda_n,
 -\lambda_1,-\lambda_2,\ldots,-\lambda_n).
\end{equation*}
so $A$ is traceless $Tr A = 0$.  Taking into account (\ref{as}) we get  
$ Tr S = 0$, i.e.  each almost K\"ahler leaf $N$ is a minimal hypersurface. 

Almost cosymplectic manifolds with K\"alerian leaves are characterized by the following theorem \cite{O2}: an almost cosymplectic manifold $(M, \varphi, \xi, \eta, g)$ has K\"ahlerian leaves if and only if
\begin{equation*}
(\nabla_X\varphi)Y = -g(\varphi AX, Y)\xi + \eta(Y)\varphi AX.
\end{equation*}

For given almost contact structure $(\varphi, \eta, \xi)$ on a  manifold $M$ we define an almost complex structure $\widetilde J$ on $M\times \mathbb R$ as
follows
\begin{equation*}
\widetilde J(X, f\dfrac{d}{dt}) = (\varphi X - f\xi, \eta(X)\dfrac{d}{dt}).
\end{equation*}
The structure $(\varphi, \xi, \eta)$ is said to be normal if $\widetilde J$ is integrable, i.e.  complex structure on $M\times\mathbb R$.  

We say that an almost cosymplectic manifold $(M,\varphi, \xi, \eta, g)$ is cosymplectic if its almost contact structure $(\varphi, \xi, \eta)$ is normal.   Cosymplectic manifolds are characterized by th condition that the tensor field $\varphi$ is parallel \cite{B}
\begin{equation*}
\nabla \varphi = 0.
\end{equation*}
Thus for cosymplectic manifold we have 
\begin{equation}
\label{commut}
R(X,Y)\varphi Z = \varphi R(X,Y)Z,
\end{equation}
that is the collineation $\varphi$ commutes with the curvature operator.   The converse statement is also true \cite{GY}, \cite{O1}:  if (\ref{commut}) is satisfied then the manifold $M$ is cosymplectic.

Now let $(N, J, G)$ be a $2n$ dimensional almost  K\"ahler manifold and $I$ a nonempty open interval.   On the product $ N\times I$ we define an almost contact metric structure $(\varphi,\xi,\eta,g)$ as follows
\begin{eqnarray*}
& \varphi(\bar X, f\dfrac{d}{dt} ) = ( J\bar X, 0), \quad\quad
\xi = \dfrac{d}{dt}, \quad\quad
\eta(\bar X, f\dfrac{d}{dt}) = f, & \\
&  g((\bar X, f\dfrac{d}{dt}),(\bar X, f\dfrac{d}{dt})) =
   G(\bar X, \bar X) + f^2, &
\end{eqnarray*}
here $f$ denotes a function on $N\times I$.  
It is simply to  verify that $(N\times I,\varphi,\xi, \eta, g)$ is an almost cosymplectic manifold and is cosymplectic if $N$ is K\"ahler.   We note that for this example we always have $A = 0$ that is the vector field $\xi$ is parallel.   

From a local point of view an almost cosymplectic manifold $M$ with vanishing tensor $A$ has a structure as described in the example above.  

\section{Conformally flat manifold}
A  Riemannian manifold $(M,g)$ is said to be locally conformally flat if  any point $p\in M$ has  a neighborhood  $U_p$ and there is a positive function $f: U_p \rightarrow U_p$ such that a metric $g'= fg$ is a flat  (Euclidean) metric on $U_p$.  

For a given Riemannian manifold $(M,g)$ the standard routine to detect the conformal flatness of $g$ is to verify that some tensor fields determined by the Riemann curvature tensor are vanishing everywhere.  If $\dim M = n \geq 4$ then one should verifies that a so called the Weyl curvature tensor vanishes.   Precisely let $R(X,Y)Z$ be the curvature operator $R(X,Y)Z = [\nabla_X, \nabla_Y]Z -\nabla_{[X,Y]}Z$ and $S$ a Ricci tensor $S(X,Y) = Tr{X \mapsto R(X,Y)Z}$.   We define a Ricci operator $Q$ requiring that $g(QX,Y) = S(X,Y)$ and a scalar curvature $s = Tr Q$.    Then the Weyl curvature $C$ is defined by   
\begin{eqnarray*}
 C(X,Y)Z & = &  R(X,Y)Z - \dfrac{1}{n-2}(g(Y,Z)QX + g(QY,Z)X \\
         & &  - g(X,Z)QY - g(QX,Z)Y) \\
         & & +\dfrac{s}{(n-1)(n-2)}(g(Y,Z)X - g(X,Z)Y).
\end{eqnarray*}
The Weyl's theorem states that the manifold $(M,g)$, $\dim M= n \geq 4$ is locally conformally if and only if the tensor $C$ vanishes.   The case of the dimension three is different for in this dimension  the Weyl curvature vanishes identically.  For $\dim M = 3$ we define a Weyl-Schouten tensor $L$
\begin{equation*}
  LX = QX - \dfrac{s}{4}X.
\end{equation*}
Then $(M,g)$ is locally conformally flat if and only if
\begin{equation}
\label{nsch}
(\nabla_XL)Y = (\nabla_YL)X.
\end{equation}

\section{Conformally flat three dimensional almost cosymplectic manifolds}
  
Let $(M, \varphi, \xi, \eta, g)$ be a three dimensional almost cosymplectic manifold.  Near a point $p\in M$ we fix an orthonormal frame $(E_1, E_2, E_3)$ of vector fields
\begin{equation*}
\begin{array}{c}
     E_1 = \xi, \quad\quad \varphi E_2 = E_3,
     \quad \quad \varphi E_3 = - E_2,   \\ 
     AE_2 = -\lambda E_2, \quad\quad AE_3 = \lambda E_3.  
\end{array}
\end{equation*}
Note that locally such frame always exists.  Moreover if $\lambda \neq 0$ it is determined uniquely up to the change of sign $(E_1, E_2) \mapsto (-E_1, -E_2)$.   The conditions $d\eta = d\varPhi = 0$ imply that the commutators $[E_i,E_j]$ should satisfy the following relations \cite{O2}
\begin{equation}
\label{comm}
\begin{array}{c}
    [ E_1, E_2 ]    =  -\lambda E_2 + \alpha E_3, \quad\quad
     [ E_1, E_3 ]  =  \alpha E_2 - \lambda E_3 \\ [+4pt]
 [ E_2, E_3 ] = \beta E_2 - \gamma E_3.
\end{array}
\end{equation}
The Jacobi identity yields the following additional conditions
\begin{equation}
\label{jac}
\begin{array}{c}
  E_2\lambda - E_3\alpha + E_1 \gamma -\alpha\beta + \gamma\lambda  =  0 ,\\
  E_3\lambda - E_2\alpha - E_1 \beta -\alpha\gamma + \beta\lambda  =  0.
\end{array}
\end{equation} 
Note that these system posses an interesting symmetry properties.  However the detailed discussion   is out of the scope of this paper.

With respect to this frame we obtain the components $S_{ij}$ of the Ricci tensor \cite{O2}
\begin{equation*}
\begin{array}{l}
 S_{11}  =  -2\lambda^2, \quad S_{12} = E_2\lambda+2\gamma\lambda, \quad
 S_{13} = -(E_3\lambda+2\beta\lambda), \\
 S_{22}  =  -E_1\lambda -E_2\gamma -E_3\beta -\beta^2 - \gamma^2, \quad
 S_{23}  = -2\alpha\lambda, \\
\nonumber S_{33}  =  E_1\lambda- E_2\gamma - E_3\beta -\beta^2 -\gamma^2.
\end{array}
\end{equation*}
and the scalar curvature $s= S_{11}+S_{22}+S_{33}$
\begin{equation*}
  s =  -2E_2\gamma -2E_3\beta -2(\beta^2 +\gamma^2+\lambda^2).
\end{equation*}

Now let assume  that $M \subset \mathbb R^3$ is a domain.  Let $p\in M$, $p =(x,y,z)$ and
\begin{equation*}
 E_1 = \xi = \dfrac{\partial}{\partial z}, \quad 
  E_2 = (a^1,a^2, a^3),\quad 
 E_3 = (b^1, b^2, b^3), 
\end{equation*}
where $a^i, b^i$ are some functions on $M$ and 
\begin{equation*}
  (c^1,c^2,c^3) = c^1 \dfrac{\partial}{\partial x} +
                  c^2 \dfrac{\partial}{\partial y} +
                  c^3 \dfrac{\partial}{\partial z}
\end{equation*}
Now (\ref{comm}), (\ref{jac}) and (\ref{nsch}) form a nonlinear overdetermined system of differential equations of the second order with respect to the unknown functions $a^i, b^j, \alpha, \beta, \gamma ,\lambda $.  

As we mentioned in the Introduction all cosymplectic locally conformally flat manifolds  can be described completely.   Indeed we known that a manifold $M$ of this type is locally a product of an nonempty open interval and the a two dimensional K\"ahler manifold $N$.  Therefore $N$ is of constant sectional curvature.  

A  non-cosymplectic case is more complicated.   Note that the example described in the Introduction can be slightly generalized.  Let $U \subset \mathbb R^3$ be a domain in $\mathbb R^3 = \{ (x,y,z) | x,y,z \in \mathbb R\}$.  Let $f = f(x,z) > 0$, $u = u(x) > 0$ be a real functions on $U$. We define an almost contact metric structure $(\varphi, \xi, \eta, g)$ as follows
\begin{eqnarray*}
\nonumber 
  & g = f(x,z)^2 dx^2 + \dfrac{u(x)^2}{f(x,z)^2}dy^2 + dz^2, & \\
  & \xi = \dfrac{\partial}{\partial z}, \quad \quad \eta = dz, & \\
\nonumber 
  & \varphi \dfrac{\partial}{\partial x} = 
         \dfrac{f(x,z)^2}{u(x)}\dfrac{\partial}{\partial y}, \quad\quad
    \varphi \dfrac{\partial}{\partial y} = 
    - \dfrac{u(x)}{f(x,z)^2} \dfrac{\partial}{\partial x}. &
\end{eqnarray*}
For the fundamental form we have
\begin{equation*}
\varPhi = 2 u(x)dx\wedge dy.
\end{equation*}
Obviously $d\eta=d\varPhi=0$.   
Rearranging terms we  write down $g$ as follows
\begin{equation*}
g = \dfrac{u(x)^2}{f(x,z)^2}\hspace{2pt}( dy^2 + \dfrac{f(x,z)^4}{u(x)^2}dx^2 +      \dfrac{f(x,z)^2}{u(x)^2}dz^2).
\end{equation*}
It is evident that $g$ is conformally flat if and only if the term inside the parenthesis  is  a conformally flat metric. 
\begin{proposition}
The almost cosymplectic manifold $(U,\varphi, \xi, \eta, g)$  is conformally flat if and only if the functions $f$ and $u$ satisfy the following differential equation
\begin{equation}
\label{constsec}
   2\partial_z^2 f - \partial_x^2\dfrac{1}{f} - 
   \partial_x(\dfrac{\partial_x \ln u}{f}) = -\kappa \dfrac{f^3}{u^2},
\end{equation}
for a constant $\kappa$.
\end{proposition}
\begin{proof}
The metric
\begin{equation*}
dy^2 + \dfrac{f(x,z)^4}{u(x)^2}dx^2 +  \dfrac{f(x,z)^2}{u(x)^2}dz^2,
\end{equation*}
is conformally flat if and only if 
\begin{equation*}
\dfrac{f(x,z)^4}{u(x)^2}dx^2 +  \dfrac{f(x,z)^2}{u(x)^2}dz^2,
\end{equation*}
 is of constant sectional curvature $\kappa$. The latter is equivalent to the functions $f$ and $u$ should satisfy (\ref{constsec}).
\end{proof}

  We note that the function $u$ in this equation plays a role of the functional parameter and is not obvious that for a given $u$ the  solution exists.  However if $u$ is real analytic then we can apply the Cauchy-Kovalevska theorem. 

\vspace{0.5cm}
\noindent{\bf Examples.}
Setting  $f(x,z) = \dfrac{t(z)}{s(x)}$, $t(z)> 0$, $s(x)> 0$ in (\ref{constsec}) we obtain 
\begin{equation*}
2t^{''}_z -\dfrac{s((s\cdot u)'_x/u)_x'}{t} = -\kappa \dfrac{t^3}{(s\cdot u)^2}.
\end{equation*}
 If $\kappa = 0$ then we have two independent equations
\begin{equation*}
 t^{''}_z = 0, \quad \quad ((s\cdot u)'_x/u)_x'= 0. 
\end{equation*}
Solving these equations one gets
\begin{equation}
  g = \dfrac{u(x)^2(Az+B)^2}{(C\int u(x) +D)^2}dx^2 +
      \dfrac{(C\int u(x) +D)^2}{(Az+B)^2}dy^2 + dz^2, 
\end{equation}
as we see $u$ may be arbitrary function.  Note that if $u = e^{ax}$ then for properly chosen constants $A, B, C, D$ we  obtain the example described in the Introduction. 

If $\kappa \neq 0$ a solution of the form $t(z)/s(x)$ exists  only if $s\cdot u = C = const > 0$.  In this case  
\begin{equation*}
  t^{''}_z = -\dfrac{\kappa}{2C^2} t^3,
\end{equation*}
which yields 
\begin{equation*}
z  = \pm 2C\int \dfrac{dt}{\sqrt{D -\kappa t^4}}.
\end{equation*}
The constant $D$ should be properly chosen depending on $\kappa$ and  a domain of definition of $t$.  Explicitly the metric is given
by
\begin{equation}
g =(C t)^2 dv^2 + \dfrac{1}{(Ct)^2}dy^2 +\dfrac{4C^2}{D -\kappa t^4}dt^2,
\end{equation}
where $dv = u(x)dx$. As above the function $u$ can be arbitrary. 
\vspace{0.5cm}

The examples  above  have some interesting properties.  We note which is evident that the field $\partial_y$ is Killing.  What is more important this is the field of eigenvectors
\begin{equation*}
A\dfrac{\partial}{\partial y} = 
\partial_z {\rm ln} f \dfrac{\partial}{\partial y}.
\end{equation*}
From this point of view without reference  to any particular local chart we may say that these manifolds have the property that there is a smooth field of eigenvectors of the operator $A$ which is Killing. 

\begin{theorem}
Let $(M,\varphi, \xi, \eta , g)$ be a three dimensional almost cosymplectic
manifold.  Assume that a smooth field $K$ of eigenvectors $K$ of the operator 
$A = - \nabla\xi$, $A\neq 0$ everywhere, is a Killing vector field.   Then the manifold $M$ is conformally flat if and only if the  Laplacian 
$\Delta \dfrac{1}{|K|}$ of the inverse of the length of the field $K$ satisfies the following equation
\begin{equation}
\Delta \dfrac{1}{|K|} = \dfrac{\kappa}{2} \dfrac{1}{|K|^3},
\end{equation}
for a constant $\kappa$. 
\end{theorem}  
\begin{proof}
Note that $[\xi,K]=0$.  
Indeed, let $X$ be arbitrary vector field. Then  as $K$ is Killing and 
$A$ symmetric we find 
\begin{equation}
  0 =  (\mathcal L_K g)(\xi, X) =  g([\xi,K], X).
\end{equation}
 Now we introduce a local orthonormal frame $(E_1, E_2, E_3)$ and 
\begin{equation}
\label{xikill}
E_1 = \xi, \quad\quad \varphi E_2 = E_3 ,  \quad \quad E_3  = K/|K|.
\end{equation} 
As  $A E_2 = -\lambda E_2$,  $A E_3 = \lambda E_3$ by (\ref{xikill}) and (\ref{comm}) 
\begin{equation*}
[E_1, E_2 ] =  -\xi\ln |K| E_2,  \quad\quad 
[E_1, E_3 ] = [\xi , E_3] = \xi\ln|K| E_3.
\end{equation*}
 As a collorary we get $\lambda = \xi\ln |K|$.  Therefore all the distributions $(E_i, E_j)$, $i < j$ are involutive.  Equivalently\begin{equation*} 
\theta^1\wedge d\theta^1 = \theta^2\wedge d\theta^2 =
\theta^3\wedge d\theta^3 = 0
\end{equation*}
where $\theta^i$ are dual forms.   Hence there are locally nonzero functions $u_i$ such that the frame $u_i E_i$ is holonomic, i.e $[u_i E_i , u_j E_j] = 0$.  The functions $u_i$ are simply  ``integrability factors'' $d(\theta^i/u_i) = 0$.
It is clear that as $d\theta^1 = d\eta = 0$ we may assume that $u_1 = 1$.  We introduce the following denotations
\begin{equation*}
\xi = \dfrac{\partial}{\partial z} = \partial_z, \quad\quad
u_2E_2 = \partial_x, \quad\quad u_3E_3 = \partial_y.
\end{equation*}
 The $\partial_x, \partial_y, \partial_z$ are eigenvectors fields and as they correspond to different eigenvalues, resp. $-\lambda$, $\lambda$, $0$ they are
pairwise orthogonal.  Therefore the  metric takes the following form  
\begin{equation}
\label{metr}
g  = f^2dx^2+\dfrac{u^2}{f^2}dy^2+dz^2,
\end{equation}
for a functions $f > 0 $, $u > 0$.
It is clear that  $K=\beta\partial_y$ with a nonzero coefficient $\beta$.  We always can assume that $\beta = 1$. Indeed the conditions
\begin{equation*}
 0 = [\xi , K] = [\partial_z, K], \quad \quad 0= (\mathcal L_Kg)(\partial_x, \partial_y),
\end{equation*}
and (\ref{metr}) imply  $K= \beta(y)\partial_y =\partial_{y'}$.  
For $\partial_y$ is Killing and $d\varPhi = 2d(\theta^2\wedge\theta^3) =0$ we obtain $f = f(x,z)$, $u = u(x)$. 

 Let $\Delta$ denote the metric  Laplacian, i.e.
\begin{equation*}
  \Delta v = - Tr \{X \mapsto \nabla_X \mbox{grad} v\},
\end{equation*}
where $dv(X) = g(X, \mbox{grad} v)$.  
 Taking into account $u=u(x)$ one verifies that 
\begin{equation*}
-u \Delta \dfrac{f}{u} =  \partial^2_z -\partial^2_x\dfrac{1}{f}
- \partial_x (\dfrac{\partial_x \ln u}{f}).
\end{equation*}
Therefore by the above identity the condition (\ref{constsec}) can be written as\begin{equation*}
 \Delta\dfrac{f}{u} = \dfrac{\kappa}{2} (\dfrac{f}{u})^3.
 \end{equation*}
Finally  $f/u = 1/|K|$.
\end{proof}

\end{document}